\documentclass{article}

\usepackage{PRIMEarxiv}

\usepackage[utf8]{inputenc} 
\usepackage[T1]{fontenc}    
\usepackage{hyperref}       
\usepackage{url}            
\usepackage{booktabs}       
\usepackage{amsfonts}       
\usepackage{nicefrac}       
\usepackage{microtype}      
\usepackage{lipsum}
\usepackage{fancyhdr}       
\usepackage{graphicx}       
\graphicspath{{media/}}     
\usepackage{amsmath}
\usepackage{amsthm}
\usepackage{fdsymbol}
\usepackage{tikz}  
\usepackage{float}
\usepackage{chemfig}
\usepackage[version=3]{mhchem}

\theoremstyle{plain} 
\newtheorem{theorem}{Theorem}[section]
\theoremstyle{definition}
\newtheorem{corollary}[theorem]{Corollary}

\newtheorem{proposition}{Proposition}[section]
\newtheorem{definition}{Definition}[section]
\newtheorem{remark}{Remark}[section]
\newtheorem{example}{Example}[section]

\title{A study on $k$-coalescence of two graphs
}

\author{
  Najiya V K, Chithra A V \\
 Department of Mathematics \\
 National Institute of Technology, Calicut\\
 Kerala, India-673601\\
  \texttt{najiya\_p190046ma@nitc.ac.in, chithra@nitc.ac.in} \\
}

\begin{document}
\maketitle

\begin{abstract}
    The $k$-coalescence of two graphs is obtained by merging a $k$-clique of each graph. The $A_\alpha$-matrix of a graph is the convex combination of its degree matrix and adjacency matrix. In this paper, we present some structural properties of a non-regular graph which is obtained from the $k$-coalescence of two graphs. Also, we derive the $A_\alpha$-characteristic polynomial of $k$-coalescence of two graphs and then compute the $A_\alpha$-spectra of $k$-coalescence of two complete graphs. In addition, we estimate the $A_\alpha$-energy of $k$-coalescence of two complete graphs. Furthermore, we obtain some topological indices of vertex coalescence of two graphs, and as an application, we determine some indices of some family of graphs. From these results, we calculate the Wiener index, hyper-Wiener index etc. of the organic compound 1,2-dicyclohexylethane(\ce{C_{14}H_{26}}).
\end{abstract}

\textbf{AMS classification}: 05C50
\newline
\\
{\bf{Keywords}}: {\it{ $A_\alpha$-spectrum, Coalescence, Wiener index, hyper-Wiener index, Zagreb index}}

\section{Introduction}\label{intro}

Let $G$ be a simple graph on $n$ vertices with vertex set ${v_1,v_2,\dots,v_n}$ and $m$ edges. The adjacency matrix\cite{chartrand2006introduction} $A(G)=[a_{ij}]$ of $G$ is defined as an $n\times n$ matrix with $a_{ij}=1$ if $v_i$ and $v_j$ are adjacent, $0$ otherwise. The signless Laplacian matrix $Q(G)$ of $G$ has the form $D(G)+A(G)$, where $D(G)$ is a diagonal matrix with $a_{ii}=deg(v_i)$. In \cite{nikiforov2017merging}, Nikiforov introduced a new matrix, which is a convex combination of $D(G)$ and $A(G)$, defined as $A_\alpha(G)=\alpha D(G)+(1-\alpha)A(G)$, where $\alpha\in[0,1]$. The $A_\alpha$ matrix, $A_\alpha(G)$ coincides with $A(G)$, $D(G)$ and $\frac{1}{2}Q(G)$ when $\alpha=0,1,\frac{1}{2}$ respectively.

For a matrix $M$, $\Phi(M,\lambda)$ denotes the characteristic polynomial of $M$. The solution for this polynomial constitutes the spectrum of $M$. 
The adjacency energy $\varepsilon(G)$ of a graph $G$ is defined as the sum of absolute values of its adjacency eigenvalues. If $\lambda_i(A_\alpha(G))$ denotes the $A_\alpha$-eigenvalues of $G$, then the $A_\alpha$-energy\cite{pirzada2021alpha} is defined as $\varepsilon_\alpha=\displaystyle\sum_{i=1}^n\left|\lambda_i(A_\alpha(G))-\frac{2\alpha m}{n}\right|$. If $G$ is a regular graph then $A_\alpha$-energy is $(1-\alpha)\varepsilon(G)$.

Let $G_1$ and $G_2$ be two graphs on $n_1, n_2$ vertices and $m_1, m_2$ edges. The $k$-coalescence\cite{sudhir2021spectra} $G_1\circ_k G_2$ of $G_1$ and $G_2$ is the graph obtained by merging a clique of order $k$ of both $G_1$ and $G_2$. The graph $G_1\circ_k G_2$ is non-regular with $n_1+n_2-k$ vertices and $m_1+m_2-\frac{k(k-1)}{2}$ edges. If $k=1$, it is called the vertex coalescence and if $k=2$, it is called the edge coalescence\cite{jog2016adjacency}. The merged clique of order $k$ is represented by $\mathcal{Q}$. It is difficult to calculate a general formula for $A_\alpha$-energy of non-regular graphs. In this paper, we obtain a formula for the $A_\alpha$-energy of vertex coalescence and edge coalescence of two complete graphs.

A topological index is a real number that is invariant under graph isomorphism and is derived from the structure of a graph. They have become prevalent due to their applications in several areas, including chemistry and networks. The most famous indices are Zagreb, Randić, Wiener, harmonic indices and their variants. Many chemists and mathematicians have extensively studied the Wiener index. In this paper, we compute certain topological indices, such as the Wiener index, hyper Wiener index, etc., of $k$-coalescence of two graphs.

Throughout this paper, $K_n$ denotes the complete graph of order $n$. The matrix $I_n$ denotes the identity matrix of order $n$, $O_{m\times n}$
denotes the $0$ matrix of order $m\times n$ and $J_{m\times n}$ is the matrix
of order $m\times n$ with all entries equal to one.

This paper is organised as follows. Section \ref{prelim} presents some definitions and results used for our work. In Section \ref{stru}, we determine some structural properties of $k$-coalescence of two graphs. In Section \ref{kcoal}, we estimate the $A_\alpha$-characteristic polynomial of $k$-coalescence of two graphs. In Section \ref{kcoalcom}, $A_\alpha$-spectrum and $A_\alpha$-energy of $k$-coalescence of two complete graphs are determined. In Section \ref{top}, some topological indices of vertex coalescence of two graphs are computed.

\section{Preliminaries}\label{prelim}
This section presents some definitions and theorems used to prove the main results.

\begin{definition}\cite{chartrand2006introduction}
	The distance $d(u, v)$ between two vertices $u$ and $v$ in $G$ is the length of the shortest path joining them, if any; otherwise, $d(u, v)=\infty$.
\end{definition}

\begin{definition}\cite{chartrand2006introduction}
	The girth $g(G)$ of a graph $G$ is defined as the length of the shortest cycle (if any) in $G$.
\end{definition}	

\begin{definition}\cite{chartrand2006introduction}
	A complete subgraph of $G$ is called a clique of $G$, and a clique of $G$ is
a maximal clique of $G$ if it is not properly contained in another clique of $G$.  The clique number of a graph $G$ is the number of vertices in a maximal clique of $G$, denoted by $\omega(G)$.

\end{definition}

\begin{definition}\cite{chartrand2006introduction}
	The vertex connectivity $\mathcal{K}(G)$ of a graph $G$ is the minimum number of vertices whose removal results in a disconnected or trivial graph.
\end{definition}

\begin{definition}\cite{chartrand2006introduction}
    An edge-cut in a nontrivial graph $G$ is a set $X$ of edges of $G$
such that $G - X$ is disconnected. The edge connectivity $\lambda(G)$ of a nontrivial graph $G$ is the cardinality of a minimum edge-cut of $G$.
\end{definition}

\begin{definition}\cite{chartrand2006introduction}
	Let $G$ be a graph. A circuit $C$ of $G$ that contains every edge of $G$ is an Eulerian circuit. A connected graph $G$ is called an Eulerian if $G$ contains an Eulerian circuit.
\end{definition}

\begin{theorem}\label{euler}\textnormal{\cite{chartrand2006introduction}}
    A nontrivial connected graph $G$ is Eulerian if and only if every vertex of $G$ has an even degree.
\end{theorem}

\begin{definition}\cite{chartrand2006introduction}
	If a graph $G$ has a spanning cycle $C$, then $G$ is called a Hamiltonian graph and $C$ a Hamiltonian cycle.

 If a graph $G$ contains a cut vertex, then $G$ cannot be Hamiltonian.
\end{definition}

\begin{definition}\cite{chartrand2006introduction}
	A set of vertices in a graph is independent if no two of them are adjacent. The largest number of vertices in such a set is called the independence number of $G$, and it is denoted by $\beta_{0}(G)$ or $\beta_{0}$.
\end{definition}

\begin{definition}\cite{chartrand2006introduction}
A proper vertex colouring of a graph $G$ is an assignment of colours to the vertices of $G$, one colour to each vertex so that adjacent vertices are coloured differently. A graph $G$ is $k$-colourable if a colouring of $G$ exists from a set of $k$ colours. The minimum positive integer $k$ for which $G$ is $k$-colourable is the chromatic number of $G$ and is denoted by $\chi(G)$.
\end{definition}

\begin{definition}\label{W}\cite{dobrynin2002wiener}
	Let $G$ be a finite, undirected, connected simple graph. Wiener index $W(G)$ of a graph $G$ is a distance based topological index, defined as the sum of the distance between all pairs of vertices in a graph $G$. Let $d_{G}(v)$ be the sum of distance between $v$ and all other vertices of $G$, then $$W(G)=\displaystyle\sum_{\{u, v\}\subseteq V(G)}d(u, v)=\frac{1}{2}\sum_{v\in V(G)}d_{G}(v).$$	
\end{definition}
\begin{definition}\label{WW}\cite{khalifeh2008hyper}
	Let $G$ be a finite, undirected, connected simple graph. The hyper-Wiener index $WW(G)$ of a graph $G$ is defined as $$WW(G)=\displaystyle\frac{1}{2}W(G)+\frac{1}{2}\sum_{\{ u, v\}\subseteq V(G)}d^{2}(u, v),$$ where $d^{2}(u, v)=d(u, v)^{2}$ and $d(u,v)$ is distance from $u$ to $v$. Let $d^{2}_{G}(v)$ be the sum of square of distances between $v$ and all other vertices of $G$, then $$WW(G)=\displaystyle\frac{1}{2}W(G)+\frac{1}{4}\sum_{v\in V(G)}d^{2}_{G}(v).$$	
\end{definition}

\begin{definition}\label{F}\cite{furtula2015forgotten}
	The forgotten topological index $F(G)$ of a graph $G$ is
	$$F(G)=\displaystyle\sum_{v\in V(G)}deg(v)^{3}=\sum_{uv\in E(G)}\left(deg(u)^{2}+deg(v)^{2}\right).$$
\end{definition}

\begin{definition}\label{M1}\cite{nikolic2003zagreb}
	The first Zagreb index $M_{1}(G)$ of a graph $G$ is $M_{1}(G)=\displaystyle\sum_{v\in V(G)}deg(v)^{2}$.
\end{definition}




\begin{definition}\label{NK}\cite{narumi1984simple}
	The Narumi - Katayama index $NK(G)$ of a graph $G$ is $NK(G)=\displaystyle\prod_{v\in V(G)}deg(v)$.
\end{definition}




\section{Structural properties of \texorpdfstring{$k$}{}-coalescence of graphs}\label{stru}
This section estimates the structural properties of $k$-coalescence of graphs, namely, chromatic number, vertex connectivity, edge connectivity, etc. Throughout the section, $G_i$ represents graphs on $n_i$ vertices.

We represent a graph's maximum degree and minimum degree by $\Delta(G)$ and $\delta(G)$, respectively.
\begin{proposition}
Let $G_i$ be regular graphs of order $n_i$ and regularity $r_i$ for $i=1,2$ and let $G=G_1\circ_kG_2$. Then $\Delta(G)=r_1+r_2-k+1$. 

If $k=n_1$ or $n_2$, then $\delta(G)=max\{r_1,r_2\}$ and if $k<n_1,n_2$, then $\delta(G)=min\{r_1,r_2\}$.
\end{proposition}
\begin{proof}
Let $v$ be any vertex of $G_1\circ_kG_2$. Then
$$deg(v)=\begin{cases}
deg_{G_1}(v)& \text{if }v\in V(G_1\setminus \mathcal{Q}),\\
deg_{G_2}(v)& \text{if } v\in V(G_2\setminus \mathcal{Q}),\\
deg_{G_1}(v)+deg_{G_2}(v)-k+1 &\text{if } v\in \mathcal{Q}.
\end{cases}$$

If $G_1$ and $G_2$ are regular, then the vertices in $\mathcal{Q}$ have degree $r_1+r_2-k+1$, which is greater than $r_1$ and $r_2$. Thus the maximum degree, $\Delta(G)=r_1+r_2-k+1$.

Without loss of generality, assume that $k=n_1$ and $n_1<n_2$. Then all the vertices in $G_1$ will be merged to a $k$ clique in $G_2$ resulting in $G_2$ itself. Then $\delta(G)=r_2=max\{r_1,r_2\}$.

Next assume $k<n_1,n_2$. Then there are vertices of degrees $r_1$ and $r_2$ in $G_1\circ_kG_2$. Thus $\delta(G)=min\{r_1,r_2\}$.
\end{proof}

\begin{proposition}
Let $g(G_i)$ be the girth of $G_i,i=1,2$. Then the girth of $G_1\circ_kG_2$ $$g(G_1\circ_kG_2)=\begin{cases}
3 & \text{ if } k\geq 3,\\
min\{g(G_1),g(G_2)\} &\text{ if } k \leq 2.
\end{cases}$$
\end{proposition}
\begin{proof}
If $k$ is greater than $2$, then the graph $G_1\circ_kG_2$ will have a cycle of length $3$ in $\mathcal{Q}$.

If $k\leq2$, then the shortest cycle in $G_1\circ_kG_2$ will be the shortest cycle in either $G_1$ or $G_2$.
\end{proof}

\begin{proposition}
Let $\omega_i$ be the clique number of $G_i,i=1,2$. Then the clique number of $G_1\circ_kG_2$,
$$\omega(G_1\circ_kG_2)=max\{\omega_1,\omega_2\}.$$
\end{proposition}
\begin{proof}
The graph $G=G_1\circ_kG_2$ has $G_1$ and $G_2$ as induced subgraphs. Thus, any clique of $G_1$ and $G_2$ is a clique of $G$ as well. Also, the merging of vertices does not produce a new clique. Hence, $\omega(G_1\circ_kG_2)=max\{\omega_1,\omega_2\}.$
\end{proof}

\begin{proposition}
Let $\mathcal{K}_i$ be the vertex connectivity of $G_i,i=1,2$. Then the vertex connectivity of $G_1\circ_kG_2$,
$$\mathcal{K}(G_1\circ_kG_2)=min\{\mathcal{K}_1,\mathcal{K}_2,k\}.$$
\end{proposition}
\begin{proof}
Suppose $\mathcal{K}_1$ and $\mathcal{K}_2$ are greater than or equal to $k$, then $G_1\circ_kG_2$ can be disconnected by removing $k$ vertices in $\mathcal{Q}$. Otherwise, the minimum vertex-cut of $G_i$ belongs to $V(G_i\setminus \mathcal{Q})$. Therefore, the vertex connectivity of $G_1\circ_kG_2=min\{\mathcal{K}_1,\mathcal{K}_2,k\}$.
\end{proof}

\begin{proposition}
Let $\lambda_i$ be the edge connectivity of $G_i,i=1,2$. Then the edge connectivity of $G_1\circ_kG_2$,
$$\lambda_i(G_1\circ_kG_2)=min\{\lambda_1,\lambda_2\}.$$
\end{proposition}
\begin{proof}
If the minimum edge-cut of $G_1$ and $G_2$ does not belong to $\mathcal{Q}$, then edge connectivity of $G_1\circ_kG_2=min\{\lambda_1,\lambda_2\}$. If the minimum edge-cut of $G_1$ or $G_2$ is in $\mathcal{Q}$, then it is same as the minimum edge-cut of $G_1\circ_kG_2$, therefore $G_1\circ_kG_2=min\{\lambda_1,\lambda_2\}$.
\end{proof}

\begin{proposition}
Let $G_1$ and $G_2$ be Eulerian graphs. Then the graph $G_1\circ_kG_2$ is Eulerian if and only if $k$ is odd.
\end{proposition}
\begin{proof}
If $G_1$ and $G_2$ are Eulerian, then by Theorem \ref{euler}, every vertex of $G_1$ and $G_2$ are of even degree. For a vertex $v$ in $G_1\circ_kG_2$
$$deg(v)=\begin{cases}
deg_{G_1}(v)& \text{if }v\in V(G_1\setminus \mathcal{Q}),\\
deg_{G_2}(v)& \text{if } v\in V(G_2\setminus \mathcal{Q}),\\
deg_{G_1}(v)+deg_{G_2}(v)-k+1 &\text{if } v\in \mathcal{Q}.
\end{cases}$$
Then $G_1\circ_kG_2$ is Eulerian if and only if $deg_{G_1}(v)+deg_{G_2}(v)-k+1$ is even, that is $k$ is odd.
\end{proof}

\begin{figure}[H]
    \centering
  \begin{tikzpicture}  
  [scale=0.9,auto=center]
   \tikzset{dark/.style={circle,fill=black}}

 \node[dark] (a1) at (0,0){} ;  
  \node[dark] (a2) at (0,2)  {}; 
  \node[dark] (a3) at (2,2)  {};  
  \node[dark] (a4) at (2,0) {};

  \draw (a1) -- (a2);
  \draw (a2) -- (a3);  
  \draw (a3) -- (a4);  
  \draw (a4) -- (a1);

 \node at (1,-1) {$C_4$};

 \node[dark] (c1) at (4,0){} ;  
  \node[dark] (c2) at (6,0)  {}; 
  \node[dark] (c3) at (8,0)  {};  
  \node[dark] (c4) at (4,2) {};  
  \node[dark] (c5) at (6,2)  {};  
  \node[dark] (c6) at (8,2)  {};  

  \draw (c1) -- (c2);
  \draw (c2) -- (c3);  
  \draw (c3) -- (c6);  
  \draw (c4) -- (c5);  
  \draw (c5) -- (c6);  
  \draw (c1) -- (c4);
  \draw (c2) -- (c5);

  \node at (6,-1) {$C_4\circ_2 C_4$};

 \node[dark] (d1) at (10,1.5){} ;  
  \node[dark] (d2) at (11.5,0)  {}; 
  \node[dark] (d3) at (13,1.5)  {};  
  \node[dark] (d4) at (11.5,3) {};  
  \node[dark] (d5) at (14.5,0)  {}; 
  \node[dark] (d6) at (14.5,3) {};  
   \node[dark] (d7) at (16,1.5){} ;

   \draw (d1) -- (d2);
  \draw (d2) -- (d3);  
  \draw (d3) -- (d4);  
  \draw (d4) -- (d1);  
  \draw (d3) -- (d5);  
  \draw (d5) -- (d7);  
  \draw (d7) -- (d6);
  \draw (d3) -- (d6);

  \node at (13,-1) {$C_4\circ_1 C_4$};

\end{tikzpicture}  
    \caption{$C_4\circ_2 C_4$ is not Eulerian whereas $C_4\circ_1 C_4$ is Eulerian.}
    \label{eulerian}
\end{figure}
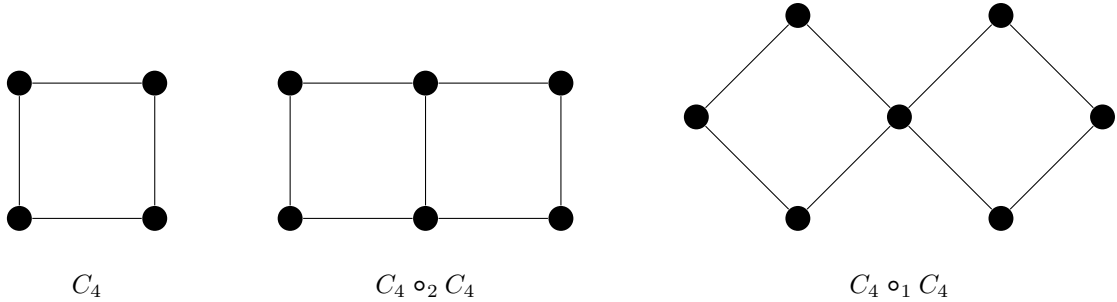

\begin{proposition}
For $k>1$, the graph $G_1\circ_kG_2$ is Hamiltonian if and only if both $G_1$ and $G_2$ are Hamiltonian. If $k=1$, then $G_1\circ_kG_2$ is not Hamiltonian.
\end{proposition}
\begin{proof}
If $k=1$, the vertex in $\mathcal{Q}$ is a  vertex cut. Then $G_1\circ_kG_2$ is not Hamiltonian.

Consider $k\geq2$. Let $n_i$ be the order of $G_i$, $i=1,2$. Assume $G_1$ and $G_2$ are Hamiltonian, then they have a Hamiltonian cycle $u_1u_2\cdots u_{n_1}u_1$ and $v_1v_2\cdots v_{n_2}v_1$ respectively, where $u_i$'s are the vertices of $G_1$ and $v_i$'s are the vertices of $G_2$. 

Let $u_r, u_{r+1}, \cdots,u_{r+k}$ and $v_1, v_{2}, \cdots, v_{k}$ be the vertices merging in $G_1\circ_kG_2$. We denote the resulting vertices as $w_1, w_{2}, \cdots,w_{k}$. The merging is in such a way that $v_1$ merge with $u_{r+m+1}$ for some $m\in\{r,r+1,\cdots,r+k\}$ and is denoted as $w_{m+1}$, $v_{2}$ merge with $u_{r+m+2}$ and is denoted as $w_{m+2}$ and so on(see Figure \ref{fig}). Then we can construct a new Hamiltonian cycle $u_1u_2\cdots w_1w_2\cdots w_m v_{k+1}v_{k+2}\cdots v_{n_{2}}w_{m+1}\cdots w_{k}\cdots u_{n_1}u_1$. Hence $G_1\circ_kG_2$ is Hamiltonian.

Conversely, if $G_1\circ_kG_2$ is Hamiltonian, then there exists a Hamiltonian cycle \\
$u_1u_2\cdots w_1w_2\cdots w_m v_{k+1}v_{k+2}\cdots v_{n_{2}}w_{m+1}\cdots w_{k}\cdots u_{n_1}u_1$.
In this cycle, consider the path $w_{m+1}\cdots w_{k}\cdots u_{n_1}u_1u_{2}\cdots w_1w_2\cdots w_{m}$. Since there is an edge between $w_m$ and $w_{m+1}$, adding this edge to the path will produce a cycle containing all the vertices of $G_1$. Therefore $G_1$ is Hamiltonian. Similarly, we can show that $G_2$ is also Hamiltonian.
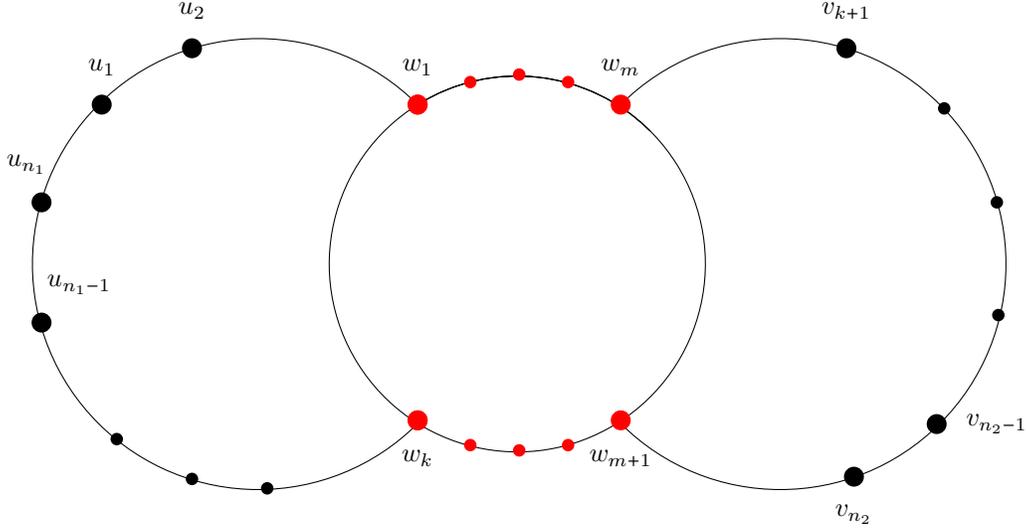
\begin{figure}[H]
    \centering
   \begin{tikzpicture}
\tikzset{dark/.style={circle,fill=black}}
\tikzset{blue/.style={circle,fill=red}}

    \draw (0,0) arc (45:315:3); 
    \node[dark,scale=0.8] (u1) at (-4.2,0) {};
    \node[dark,scale=0.8] (u2) at (-3,0.75) {};
    \node[blue,scale=0.8] (w1) at (0,0) {};
    \node[blue,scale=0.8] (wk) at (0,-4.2){};
    \node[dark,scale=0.8] (un1) at (-5,-2.9) {};
    \node[dark,scale=0.8] (un) at (-5,-1.3) {};
     \node[dark,scale=0.5] at (-2,-5.105){};
     \node[dark,scale=0.5] at (-3,-4.98){};
     \node[dark,scale=0.5] at (-4,-4.45){};
     
    \draw (0,0) arc (122:-315:2.5); 
     \node[blue,scale=0.5] at (0.7,0.3){};
      \node[blue,scale=0.5] at (1.35,0.4){};
       \node[blue,scale=0.5] at (2,0.3){};

        \node[blue,scale=0.5] at (0.7,-4.53){};
      \node[blue,scale=0.5] at (1.35,-4.6){};
       \node[blue,scale=0.5] at (2,-4.53){};
       
    \draw (2.7,0) arc (135:-135:3); 
    \node[blue,scale=0.8] (wm) at (2.7,0) {};
    \node[dark,scale=0.8] (vk1) at (5.7,0.75) {};
    \node[dark,scale=0.5] at (7,-0.05){};
     \node[dark,scale=0.5] at (7.7,-1.3){};
     \node[dark,scale=0.5] at (7.72,-2.8){};
    \node[dark,scale=0.8] (vn21)at (6.9,-4.25) {};
     \node[dark,scale=0.8] (vn2) at (5.8,-4.95) {};
      \node[blue,scale=0.8] (wm1) at (2.7,-4.2) {};
      \node[yshift=0.5cm] at (u1) {$u_1$};
      \node[yshift=0.5cm] at (u2) {$u_2$};
      \node[yshift=0.5cm] at (w1) {$w_1$};
      \node[yshift=-0.5cm] at (wk) {$w_k$};
      \node[yshift=0.5cm, xshift=-0.2cm] at (un) {$u_{n_1}$};
      \node[yshift=0.5cm, xshift=0.5cm] at (un1) {$u_{n_1-1}$};

       \node[yshift=0.5cm] at (wm) {$w_{m}$};
      \node[yshift=-0.5cm] at (wm1) {$w_{m+1}$};
      \node[yshift=0.5cm] at (vk1) {$v_{k+1}$};
      \node[xshift=0.8cm] at (vn21) {$v_{n_2-1}$};
      \node[yshift=-0.5cm] at (vn2) {$v_{n_2}$};
   
     \node[blue,scale=0.8] (w1) at (0,0) {};
    \node[blue,scale=0.8] (wk) at (0,-4.2){};

\end{tikzpicture}
    \caption{Hamiltonian cycle in $G_1\circ_kG_2$.}
    \label{fig}
\end{figure}
\end{proof}

The following proposition gives us a lower and upper bound for the independence number of $k$-coalescence of two graphs.
\begin{proposition}
Let $\beta_0(G_i)$ be the independence number of $G_i,i=1,2$. Then the independence number of $G=G_1\circ_kG_2$ satisfies $$\beta_0(G_1)+\beta_0(G_2)-2\leq\beta_0(G)\leq \beta_0(G_1)+\beta_0(G_2).$$
\end{proposition}
\begin{proof}
Let $G=G_1\circ_kG_2$\\
Case 1: Both $G_1$ and $G_2$ are complete graphs.\\
Then the vertices in $V(G_1\setminus\mathcal{Q})$ are not adjacent to vertices in $V(G_2\setminus\mathcal{Q})$. Thus $\beta_0(G)=2=\beta_0(G_1)+\beta_0(G_2)$.

Case 2: Either $G_1$ or $G_2$ is complete.\\
Without loss of generality, assume that $G_1$ is complete and $G_2$ is not. If the independent set of $G_2$ contains a vertex in $\mathcal{Q}$, then $\beta_0(G)=\beta_0(G_2)=\beta_0(G_2)+\beta_0(G_1)-1$.  If the independent set of $G_2$ does not contain a vertex in $\mathcal{Q}$, then the independent set of $G$ contains independent vertices of $G_2$ along with a vertex from $G_1\setminus\mathcal{Q}$. Thus $\beta_0(G_1\circ_kG_2)=\beta_0(G_2)+1=\beta_0(G_2)+\beta_0(G_1)$.

Case 3: Neither $G_1$ nor $G_2$ is complete.\\
If both $G_i$'s have an independent set disjoint from $\mathcal{Q}$ then their union gives the independent set for $G$, that is, $\beta_0(G)=\beta_0(G_1)+\beta_0(G_2)$.
If one of the $G_i$'s has a vertex common in its independent set and $\mathcal{Q}$, then $\beta_0(G)=\beta_0(G_1)+\beta_0(G_2)-1$.
If both $G_i$'s has vertices common in their independent set and $\mathcal{Q}$, then $\beta_0(G)=\beta_0(G_1)+\beta_0(G_2)-2$.
\end{proof}

\begin{proposition}
Let $\chi_i$ be the chromatic number of $G_i,i=1,2$. Then the chromatic number of $G_1\circ_kG_2$,
$$\chi(G_1\circ_kG_2)=k+max\{\chi_1-k,\chi_2-k\}.$$
\end{proposition}
\begin{proof}
We need $k$ different colours to colour the vertices in $\mathcal{Q}$. Since the vertices in $V(G_1\setminus \mathcal{Q})$ and $V(G_2\setminus \mathcal{Q})$ are not adjacent, they can be coloured using $max\{\chi_1-k,\chi_2-k\}$ colours. Thus chromatic number of $G_1\circ_kG_2=k+max\{\chi_1-k,\chi_2-k\}$.
\end{proof}

\section{\texorpdfstring{$A_\alpha$}{}-characteristic polynomial of \texorpdfstring{$k$}{}-coalescence of graphs}\label{kcoal}

This section computes the $A_\alpha$-characteristic polynomial of $k$-coalescence of two graphs. Using that, the $A_\alpha$-characteristic polynomial of Lollipop graphs is estimated.

Let $G$ be a graph containing a $k$-clique. Then we partition the adjacency matrix of $G$ into the form $A(G)=\begin{bmatrix}
    B & C^T\\
    C & A(K_k)
\end{bmatrix}$. 

\begin{proposition}
Let $G_1$ and $G_2$ be two graphs of order $n_1$ and $n_2$ respectively such that $n_1+n_2>3k$. 
Then the $A_\alpha$-characteristic polynomial of $G_1\circ_k G_2$ is
\begin{align*}
\Phi(A_\alpha(G_1&\circ_k G_2),\lambda)=\Phi(A_\alpha(G_1),\lambda)\Phi(A_\alpha(G_2\backslash \mathcal{Q}),\lambda)+\Phi(A_\alpha(G_2),\lambda)\Phi(A_\alpha(G_1\backslash \mathcal{Q}),\lambda)\\
&-\Phi(A_\alpha(G_1\backslash \mathcal{Q}),\lambda)\Phi(A_\alpha(G_2\backslash \mathcal{Q}),\lambda)\Big(\alpha|D_1(\mathcal{Q})-(k-1)I|+\alpha|D_2(\mathcal{Q})-(k-1)I|\\
&\hspace{3cm}+|\lambda-\alpha(D_1(\mathcal{Q})+D_2(\mathcal{Q})-(k-1)I)-(1-\alpha)A(K_k)|\Big),
\end{align*}
where $D_i(\mathcal{Q})$ represents the degree matrix of the $k$ vertices in $\mathcal{Q}$ of $G_i$, $i=1,2$.
\end{proposition}

\begin{proof}
The $A_\alpha$-matrix of $G_1\circ_k G_2$ with proper labelling has the form 
\begin{align*}
    A_\alpha(G_1\circ_k G_2)=&\begin{bmatrix}
D & R_1^T & R_2^T\\
R_1 & A_\alpha(G_1\backslash \mathcal{Q}) & O\\
R_2 & O & A_\alpha(G_2\backslash \mathcal{Q})
\end{bmatrix},\\
\intertext{where $D=\alpha(D_1(\mathcal{Q})+D_2(\mathcal{Q})-(k-1)I)+(1-\alpha)A(K_k)$ and $R_i=(1-\alpha)C_i$, where $C_i$ is the block matrix in the adjacency matrix of $G_i$. Then, }
\Phi(A_\alpha(G_1\circ_k G_2),\lambda)=&|\lambda-A_\alpha(G_1\circ_k G_2)|\\
=&\begin{vmatrix}
\lambda-D & -R_1^T & -R_2^T\\
-R_1 & \lambda-A_\alpha(G_1\backslash \mathcal{Q}) & O\\
-R_2 & O & \lambda-A_\alpha(G_2\backslash \mathcal{Q})
\end{vmatrix}\\
&\hspace{-5cm}=\begin{vmatrix}
\lambda-D & -R_1^T & -R_2^T\\
-R_1 & O & O\\
-R_2 & O & O
\end{vmatrix} + \begin{vmatrix}
\lambda-D & -R_1^T & O\\
-R_1 & O & O\\
-R_2 & O & \lambda-A_\alpha(G_2\backslash \mathcal{Q})
\end{vmatrix} +\begin{vmatrix}
\lambda-D & O & -R_2^T\\
-R_1 & \lambda-A_\alpha(G_1\backslash \mathcal{Q}) & O\\
-R_2 & O & \lambda-A_\alpha(G_2\backslash \mathcal{Q})
\end{vmatrix}.\\
\intertext{Adding and subtracting $\begin{vmatrix}
\lambda-D & O & O\\
-R_1 & \lambda-A_\alpha(G_1\backslash \mathcal{Q}) & O\\
-R_2 & O & \lambda-A_\alpha(G_2\backslash \mathcal{Q})
\end{vmatrix}$ to  $\Phi(A_\alpha(G_1\circ_k G_2),\lambda) $, we get}
\Phi(A_\alpha(G_1\circ_k G_2),\lambda)=&\begin{vmatrix}
\lambda-D & -R_1^T & O\\
-R_1 & \lambda-A_\alpha(G_1\backslash \mathcal{Q}) & O\\
-R_2 & O & \lambda-A_\alpha(G_2\backslash \mathcal{Q})
\end{vmatrix} \\
&\hspace{-2cm}- \begin{vmatrix}
\lambda-D & O & -R_2^T\\
-R_1 & \lambda-A_\alpha(G_1\backslash \mathcal{Q}) & O\\
-R_2 & O & \lambda-A_\alpha(G_2\backslash \mathcal{Q})
\end{vmatrix} +\begin{vmatrix}
\lambda-D & O & O\\
-R_1 & \lambda-A_\alpha(G_1\backslash \mathcal{Q}) & O\\
-R_2 & O & \lambda-A_\alpha(G_2\backslash \mathcal{Q})
\end{vmatrix}\\
=& \begin{vmatrix}\lambda-A_\alpha(G_2\backslash \mathcal{Q})\end{vmatrix}\begin{vmatrix}\lambda-D & -R_1^T \\
-R_1 & \lambda-A_\alpha(G_1\backslash \mathcal{Q})\end{vmatrix}\\
&\hspace{-3cm}+\begin{vmatrix}\lambda-A_\alpha(G_1\backslash \mathcal{Q})\end{vmatrix}\begin{vmatrix}\lambda-D & -R_2^T \\
-R_2 & \lambda-A_\alpha(G_2\backslash \mathcal{Q})\end{vmatrix}-\begin{vmatrix}
\lambda-D
\end{vmatrix} \begin{vmatrix}
\lambda-A_\alpha(G_1\backslash \mathcal{Q})
\end{vmatrix} \begin{vmatrix}
\lambda-A_\alpha(G_2\backslash \mathcal{Q})
\end{vmatrix}
\end{align*}

Here,
\begin{align*}
   \begin{vmatrix}\lambda-D & -R_1^T \\
-R_1 & \lambda-A_\alpha(G_1\backslash \mathcal{Q})\end{vmatrix}=&\begin{vmatrix}\lambda-\alpha(D_1(\mathcal{Q})+D_2(\mathcal{Q})-(k-1)I)-(1-\alpha)A(K_k) & -R_1^T \\
-R_1 & \lambda-A_\alpha(G_1\backslash \mathcal{Q})\end{vmatrix}\\
=&\begin{vmatrix}\lambda-\alpha D_1(\mathcal{Q})-(1-\alpha)A(K_k) & -R_1^T \\
-R_1 & \lambda-A_\alpha(G_1\backslash \mathcal{Q})\end{vmatrix}\\
&+\begin{vmatrix}-\alpha(D_2(\mathcal{Q})-(k-1)I) & -R_1^T \\
O & \lambda-A_\alpha(G_1\backslash \mathcal{Q})\end{vmatrix}\\
=&\begin{vmatrix}
\lambda-A_\alpha(G_1)
\end{vmatrix}-\alpha
\begin{vmatrix}
D_2(\mathcal{Q})-(k-1)I
\end{vmatrix}
\begin{vmatrix}
\lambda-A_\alpha(G_1\backslash \mathcal{Q})
\end{vmatrix}.
\intertext{Similarly}
\begin{vmatrix}\lambda-D & -R_2^T \\
-R_2 & \lambda-A_\alpha(G_2\backslash \mathcal{Q})\end{vmatrix}
=&\begin{vmatrix}
\lambda-A_\alpha(G_2)
\end{vmatrix}-\alpha
\begin{vmatrix}
D_1(\mathcal{Q})-(k-1)I
\end{vmatrix}
\begin{vmatrix}
\lambda-A_\alpha(G_2\backslash \mathcal{Q})
\end{vmatrix}.
\end{align*}
Therefore,
\begin{align*}
    \Phi(A_\alpha(G_1\circ_k G_2),\lambda)=&\begin{vmatrix}\lambda-A_\alpha(G_1)\end{vmatrix}\begin{vmatrix}\lambda-A_\alpha(G_2\backslash \mathcal{Q})\end{vmatrix}+\begin{vmatrix}\lambda-A_\alpha(G_2)\end{vmatrix}\begin{vmatrix}\lambda-A_\alpha(G_1\backslash \mathcal{Q})\end{vmatrix}\\
&\hspace{-4cm}-\begin{vmatrix}
\lambda-A_\alpha(G_1\backslash \mathcal{Q})
\end{vmatrix} \begin{vmatrix}
\lambda-A_\alpha(G_2\backslash \mathcal{Q})
\end{vmatrix}\Big(\alpha\left(\begin{vmatrix}
D_1(\mathcal{Q})-(k-1)I
\end{vmatrix}+\begin{vmatrix}
D_2(\mathcal{Q})-(k-1)I
\end{vmatrix}\right)+\begin{vmatrix}
\lambda-D
\end{vmatrix}\Big)\\
=&\Phi(A_\alpha(G_1),\lambda)\Phi(A_\alpha(G_2\backslash \mathcal{Q}),\lambda)+\Phi(A_\alpha(G_2),\lambda)\Phi(A_\alpha(G_1\backslash \mathcal{Q}),\lambda)\\
&\hspace{-4cm}-\Phi(A_\alpha(G_1\backslash \mathcal{Q}),\lambda)\Phi(A_\alpha(G_2\backslash \mathcal{Q}),\lambda)\Big(\alpha\left(\begin{vmatrix}
D_1(\mathcal{Q})-(k-1)I
\end{vmatrix}+\begin{vmatrix}
D_2(\mathcal{Q})-(k-1)I
\end{vmatrix}\right)+\begin{vmatrix}
\lambda-D
\end{vmatrix}\Big).
\end{align*}

\end{proof}

\begin{corollary}
Let $G_1$ and $G_2$ be two graphs of order $n_1$ and $n_2$ respectively such that $n_1+n_2>3k$. 
Then the adjacency characteristic polynomial of $G_1\circ_k G_2$ is
\begin{align*}
\Phi(A(G_1\circ_k G_2),\lambda)=&\Phi(A(G_1),\lambda)\Phi(A(G_2\backslash \mathcal{Q}),\lambda)+\Phi(A(G_2),\lambda)\Phi(A(G_1\backslash \mathcal{Q}),\lambda)\\
&-(\lambda-k+1)(x+1)^{k-1}\Phi(A(G_1\backslash \mathcal{Q}),\lambda)\Phi(A(G_2\backslash \mathcal{Q}),\lambda).
\end{align*}
\end{corollary}

\begin{remark}\label{lollipop}
The Lollipop graph, $L(m,n-1)$ is obtained from the  coalescence of a vertex from a cycle $C_m$ and a pendant vertex from a path $P_n$. The $A_\alpha$-characteristic polynomial of the Lollipop graph is, 
$\Phi(A_\alpha(L(m,n-1)),\lambda)=$
$$\Phi(A_\alpha(P_n),\lambda)\Phi(A_\alpha(P_{m-1}),\lambda)+\Phi(A_\alpha(C_m),\lambda)\Phi(A_\alpha(P_{n-1}),\lambda)-\lambda \Phi(A_\alpha(P_{n-1}),\lambda)\Phi(A_\alpha(P_{m-1}),\lambda).$$
\end{remark}
Using the Remark \ref{lollipop}, we can calculate the $A_\alpha$-characteristic polynomial of Lollipop graphs and hence find their spectrum.

\begin{figure}[H]
    \centering
  \begin{tikzpicture}  
  [scale=0.9,auto=center]
   \tikzset{dark/.style={circle,fill=black}}

 \node[dark] (d1) at (10,1.5){} ;  
  \node[dark] (d2) at (11.5,0)  {}; 
  \node[dark] (d3) at (13,1.5)  {};  
  \node[dark] (d4) at (11.5,3) {};  
   \node[dark] (d5) at (14.5,1.5)  {}; 
   \node[dark] (d6) at (16,1.5)  {}; 
   \node[dark] (d7) at (17.5,1.5)  {};

   \draw (d1) -- (d2);
  \draw (d2) -- (d3);  
  \draw (d3) -- (d4);  
  \draw (d4) -- (d1);  
  \draw (d3) -- (d5);
  \draw (d5) -- (d6);
  \draw (d6) -- (d7);


\end{tikzpicture}  
    \caption{$L(4,3)$}
    \label{lollipopfig}
\end{figure}

\begin{example}
The $A_\alpha$-characteristic polynomial of the Lollipop graph $L(4,3)$ is, 
{\small
$\Phi(A_\alpha(L(4,3)),\lambda)=\Phi(A_\alpha(P_3),\lambda)\Phi(A_\alpha(P_{3}),\lambda)+\Phi(A_\alpha(C_4),\lambda)\Phi(A_\alpha(P_{2}),\lambda)-\lambda \Phi(A_\alpha(P_{2}),\lambda)\Phi(A_\alpha(P_{3}),\lambda).$}
\end{example}

\section{\texorpdfstring{$A_\alpha$}{}-spectrum of \texorpdfstring{$k$}{}-coalescence of complete graphs}\label{kcoalcom}
In this section, we compute the $A_\alpha$-spectrum and $A_\alpha$-energy of $K_m\circ_k K_n$.
\begin{proposition}
For $m,n> 1$, the $A_\alpha$-characteristic polynomial of $K_m\circ_k K_n$ is \\
$\Phi(A_\alpha(K_m\circ_k K_n),\lambda)=\Big(\lambda-\alpha(m+n-k)+1\Big)^{k-1}\Big(\lambda-\alpha m+1\Big)^{m-k-1}\Big(\lambda-\alpha n+1\Big)^{n-k-1} \bigg(\Big(\lambda-m+1+(1-\alpha)k\Big)\Big(\lambda-n+1+(1-\alpha)k\Big)\Big(\lambda-\alpha(m+n-2k)+1-k\Big)-(1-\alpha)^2k\Big((m+n-2k)\lambda-(m+n-2k)\alpha k-(m-k)(n-k-1)-(n-k)(m-k-1)\Big)\bigg)$.
\end{proposition}

 \begin{proof}
The degree matrix of $K_m\circ_k K_n$ with proper labelling has the form
$$D(K_m\circ_k K_n)=\begin{bmatrix}
(m+n-1-k)I_{k} & O_{k\times m-k}& O_{k\times n-k}\\
O_{m-k\times k} & (m-1)I_{m-k} & O_{m-k\times n-k}\\
O_{n-k\times k} & O_{n-k\times m-k} & (n-1)I_{n-k}
\end{bmatrix}.$$
The adjacency matrix of $K_m\circ_k K_n$ has the form 
$$A(K_m\circ_k K_n)=\begin{bmatrix}
A(K_{k}) & J_{k\times m-k}& J_{k\times n-k}\\
J_{m-k\times k} & A(K_{m-k}) & O_{m-k\times n-k}\\
J_{n-k\times k} & O_{n-k\times m-k} & A(K_{n-k})
\end{bmatrix}.$$
Thus the $A_\alpha$-matrix of $K_m\circ_k K_n$ is
$A_\alpha(K_m\circ_k K_n)=$
$$\begin{bmatrix}
\beta_1 & (1-\alpha)J_{k\times m-k}& (1-\alpha)J_{k\times n-k}\\
(1-\alpha)J_{m-k\times k1} &\beta_2 & O_{m-k\times n-k}\\
(1-\alpha)J_{n-k\times k} & O_{n-k\times m-k} & \beta_3
\end{bmatrix},$$
where $\beta_1=\alpha(m+n-1-k)I_k+(1-\alpha)A(K_{k})$, $\beta_2=\alpha(m-1)I+(1-\alpha) A(K_{m-k})$ and $\beta_3=\alpha(n-1)I+(1-\alpha)A(K_{n-k})$.

Then the characteristic polynomial of $K_m\circ_k K_n$ is 

$|\lambda I-A_\alpha(K_m\circ_k K_n)|=
\begin{vmatrix}
\lambda I_k-\beta_1 & -(1-\alpha)J_{k\times m-k} & -(1-\alpha)J_{k\times n-k}\\
-(1-\alpha)J_{m-k\times k} &\lambda I_k-\beta_2 & O\\
-(1-\alpha)J_{n-k\times k} & O & \lambda I_k-\beta_3
\end{vmatrix}.$

In the above determinant, performing $C_l\xrightarrow{} C_l+\displaystyle\frac{1-\alpha}{\lambda-m+1+(1-\alpha)k}\displaystyle\sum_{i=k+1}^{m}C_i+\frac{1-\alpha}{\lambda-n+1+(1-\alpha)k}\displaystyle\sum_{j=m+1}^{m+n-k}C_j$\\ for $l=1,2,\cdots,k$ columns we get,

$|\lambda I-A_\alpha(K_m\circ_k K_n)|=\begin{vmatrix}
\beta_4 &  -(1-\alpha)J_{k\times m-k}& -(1-\alpha)J_{k\times n-k}\\
O & \lambda I-\beta_2 & O\\
O & O & \lambda I-\beta_3
\end{vmatrix}$,\\
where $\beta_4=(\lambda-\alpha(m+n-k)+1)I_k-(1-\alpha)\left[\frac{(1-\alpha)(m-k)}{\lambda-m+1+(1-\alpha)k}+\frac{(1-\alpha)(n-k)}{\lambda-n+1+(1-\alpha)k}+1\right]J_k$.

$|\lambda I-A_\alpha(K_m\circ_k K_n)|=$
\begin{multline*}
    \left|(\lambda-\alpha(m+n-k)+1)I_k-(1-\alpha)XJ_k\right|\\
    \left|(\lambda-\alpha(m-1))I-(1-\alpha)A(K_{m-k})\right||(\lambda-\alpha(n-1))I-(1-\alpha)A(K_{n-k})|,
\end{multline*}
where $X=\left[\frac{(1-\alpha)(m-k)}{\lambda-m+1+(1-\alpha)k}+\frac{(1-\alpha)(n-k)}{\lambda-n+1+(1-\alpha)k}+1\right]$\\
Thus\\
$\Phi(A_\alpha(K_m\circ_k K_n),\lambda)=\Big(\lambda-\alpha(m+n-k)+1\Big)^{k-1}\Big(\lambda-\alpha m+1\Big)^{m-k-1}\Big(\lambda-\alpha n+1\Big)^{n-k-1} \bigg(\Big(\lambda-m+1+(1-\alpha)k\Big)\Big(\lambda-n+1+(1-\alpha)k\Big)\Big(\lambda-\alpha(m+n-2k)+1-k\Big)-(1-\alpha)^2k\Big((m+n-2k)\lambda-(m+n-2k)\alpha k-(m-k)(n-k-1)-(n-k)(m-k-1)\Big)\bigg).$

 \end{proof}

Now, in the following corollary, we obtain the $A_\alpha$-eigenvalues of $K_m\circ_k K_n$.
\begin{corollary}
 The $A_\alpha$-eigenvalues of $K_m\circ_k K_n$ are 
 \begin{enumerate}
 \item $\alpha(m+n-k)-1$ repeated $k-1$ times,
     \item $\alpha m-1$ repeated $m-k-1$ times,
     \item $\alpha n-1$ repeated $n-k-1$ times,
     \item three roots of the equation $\bigg(\Big(\lambda-m+1+(1-\alpha)k\Big)\Big(\lambda-n+1+(1-\alpha)k\Big)\Big(\lambda-\alpha(m+n-2k)+1-k\Big)-(1-\alpha)^2k\Big((m+n-2k)\lambda-(m+n-2k)\alpha k-(m-k)(n-k-1)-(n-k)(m-k-1)\Big)\bigg)=0.$
 \end{enumerate}
\end{corollary}

The following corollary helps us to determine the $A_\alpha$-energy of non-regular graph $K_m\circ_k K_n$.
\begin{corollary}
The $A_\alpha$-energy of $K_m\circ_k K_n$ is  \\
$\varepsilon_\alpha(K_m\circ_k K_n)=(k-1)\left|\alpha(1-2k)+\displaystyle\frac{2\alpha mn}{m+n-1}-1\right|+(m-k-1)\left|\alpha(1-k)+\displaystyle\frac{\alpha n(m-n+k)}{m+n-k}-1\right|+(n-k-1)\left|\alpha(1-k)+\displaystyle\frac{\alpha m(n-m+k)}{m+n-k}-1\right|+\left|\beta-X_1\right|+\left|\gamma-X_1\right|+\left|\delta-X_1\right|$,\\
where $X_1=\displaystyle\frac{\alpha\left(m^2+n^2-k^2-(m+n-k)\right)}{m+n-k}$ and $\beta, \gamma, \delta$ are roots of the equation $\bigg(\Big(\lambda-m+1+(1-\alpha)k\Big)\Big(\lambda-n+1+(1-\alpha)k\Big)\Big(\lambda-\alpha(m+n-2k)+1-k\Big)-(1-\alpha)^2k\Big((m+n-2k)\lambda-(m+n-2k)\alpha k-(m-k)(n-k-1)-(n-k)(m-k-1)\Big)\bigg)=0.$
\end{corollary}

\begin{corollary}
The $A_\alpha$-energy of $K_m\circ_k K_m$ is  \\
$\varepsilon_\alpha(K_m\circ_k K_m)=(k-1)\left|\alpha(1-2k)+\displaystyle\frac{2\alpha m^2}{2m-1}-1\right|+2(m-k-1)\left|\alpha(1-k)+\displaystyle\frac{\alpha mk}{2m-k}-1\right|+\bigg|\beta-X_2\bigg|+\bigg|\gamma-X_2|+\bigg|\delta-X_2|$,\\
where $X_2=\displaystyle\frac{\alpha\left(2m^2-k^2-(2m-k)\right)}{2m-k}$ and $\beta $, $ \gamma$ and $\delta$ are roots of the equation $(\lambda-m+1+(1-\alpha)k)^2(\lambda-2\alpha(m-k)+1-k)-(1-\alpha)^2k((2m-k)\lambda-2\alpha k(m-k)-2(m-k)(m-k-1))=0.$
\end{corollary}

\begin{corollary}
For $m,n> 1$, the $A_\alpha$-characteristic polynomial of $K_m\circ_1 K_n$ is \\
$\Phi(A_\alpha(K_m\circ_1 K_n),\lambda)=(\lambda-\alpha m+1)^{m-2}(\lambda-\alpha n+1)^{n-2} \bigg((\lambda-m+2-\alpha)(\lambda-n+2-\alpha)(\lambda-\alpha(m+n-2))-(1-\alpha)^2\Big((m+n-2)\lambda-(m+n-2)\alpha-(m-1)(n-2)-(m-2)(n-1)\Big)\bigg)$.
\end{corollary}

\begin{corollary}
 The $A_\alpha$-eigenvalues of $K_m\circ_1 K_n$ are 
 \begin{enumerate}
     \item $\alpha m-1$ repeated $m-2$ times,
     \item $\alpha n-1$ repeated $n-2$ times,
     \item three roots of the equation $(\lambda-m+2-\alpha)(\lambda-n+2-\alpha)(\lambda-\alpha(m+n-2))-(1-\alpha)^2\bigg[(m+n-2)\lambda-(m+n-2)\alpha-(m-1)(n-2)-(m-2)(n-1)\bigg]=0.$
 \end{enumerate}
\end{corollary}
The following corollary helps us to determine the $A_\alpha$-energy of a non-regular graph $K_m\circ_1 K_n$.
\begin{corollary}
The $A_\alpha$-energy of $K_m\circ_1 K_n$ is  \\
$\varepsilon_\alpha(K_m\circ_1 K_n)=\frac{m-2}{m+n-1}\bigg|\alpha n(m-n+1)-(m+n-1)\bigg|+\frac{n-2}{m+n-1}\bigg|\alpha m(n-m+1)-(m+n-1)\bigg|+|\beta-X_3|+|\gamma-X_3|+|\delta-X_3|$,\\
where $X_3=\frac{\alpha(m^2+n^2-m-n)}{m+n-1}$ and $\beta, \gamma, \delta$ are roots of the equation $(\lambda-m+2-\alpha)(\lambda-n+2-\alpha)(\lambda-\alpha(m+n-2))-(1-\alpha)^2\bigg[(m+n-2)\lambda-(m+n-2)\alpha-(m-1)(n-2)-(m-2)(n-1)\bigg]=0.$
\end{corollary}

\begin{corollary}
The $A_\alpha$-energy of $K_m\circ_1 K_m$ is  \\
$\varepsilon_\alpha(K_m\circ_1 K_m)=\frac{2(m-2)}{2m-1}(m(2-\alpha)-1)+\bigg|\frac{2m^2(1-\alpha)-5m+2-\alpha}{2m-1}\bigg|+\bigg|\beta-\frac{2m\alpha(m-1)}{2m-1}\bigg|+\bigg|\gamma-\frac{2m\alpha(m-1)}{2m-1}\bigg|$,\\
where $\beta $ and $ \gamma$ are roots of the equation $\lambda^2-(m-2+\alpha(2m-1))\lambda+2(m-1)(\alpha m-1)=0.$
\end{corollary}

\begin{corollary}
For $m,n> 2$, the $A_\alpha$-characteristic polynomial of $K_m\circ_2 K_n$ is \\
$\Phi(A_\alpha(K_m\circ_2 K_n),\lambda)=(\lambda-\alpha m+1)^{m-3}(\lambda-\alpha n+1)^{n-3}(\lambda-\alpha(m+n-2)+1)\bigg((\lambda-m+3-2\alpha)(\lambda-n+3-2\alpha)(\lambda-\alpha(m+n-4)-1)-2(1-\alpha)^2\Big((m+n-4)\lambda-(m+n-4)2\alpha-(m-2)(n-3)-(m-3)(n-2)\Big)\bigg).$
\end{corollary}

Now, in the following corollary, we obtain the $A_\alpha$-eigenvalues of $K_m\circ_2 K_n$.
\begin{corollary}
 The $A_\alpha$-eigenvalues of $K_m\circ_2 K_n$ are 
 \begin{enumerate}
     \item $\alpha m-1$ repeated $m-3$ times,
     \item $\alpha n-1$ repeated $n-3$ times,
     \item $\alpha(m+n-2)-1$,
     \item three roots of the equation $(\lambda-m+3-2\alpha)(\lambda-n+3-2\alpha)(\lambda-\alpha(m+n-4)-1)-2(1-\alpha)^2\Big[(m+n-4)\lambda-(m+n-4)2\alpha-(m-2)(n-3)-(m-3)(n-2)\Big]=0.$
 \end{enumerate}
\end{corollary}
The following corollary helps us to determine the $A_\alpha$-energy of a non-regular graph $K_m\circ_2 K_n$.
\begin{corollary}
The $A_\alpha$-energy of $K_m\circ_2 K_n$ is \\
$\varepsilon_\alpha(K_m\circ_2 K_n)=\frac{m-3}{m+n-2}\bigg|\alpha[(m-n)(n-1)+2]-1\bigg|+\frac{n-3}{m+n-2}\bigg|\alpha[(n-m)(m-1)+2]-1\bigg|+\frac{1}{m+n-2}\bigg|\alpha(2mn-3m-3n+6)-1\bigg|+|\beta-X_4|+|\gamma-X_4|+|\delta-X_4|$,\\
where $X_4=\frac{\alpha[m(m-1)+n(n-1)-2]}{m+n-2}$ and $\beta, \gamma, \delta$ are roots of the equation $(\lambda-m+3-2\alpha)(\lambda-n+3-2\alpha)(\lambda-\alpha(m+n-4)-1)-2(1-\alpha)^2\Big[(m+n-4)\lambda-(m+n-4)2\alpha-(m-2)(n-3)-(m-3)(n-2)\Big]=0.$
\end{corollary}

\begin{corollary}
The $A_\alpha$-energy of $K_m\circ_2 K_m$ is  \\
$\varepsilon_\alpha(K_m\circ_2 K_m)=\frac{m-3}{m-1}\bigg|2\alpha-1\bigg|+\frac{1}{2m-2}\bigg|\alpha(2m^2-6m+6)-1\bigg|+ \bigg|\frac{m^2(1-\alpha)-m(4-3\alpha)-\alpha+3}{m-1}\bigg|+\bigg|\beta-\frac{\alpha(m(m-1)-1)}{m-1}\bigg|+\bigg|\gamma-\frac{\alpha(m(m-1)-1)}{m-1}\bigg|$,\\
where $\beta$ and $\gamma$ are roots of the equation $\lambda^2-(m-2+2\alpha(m-1))\lambda+2\alpha m^2-m(2\alpha+3)-2\alpha+5=0.$
\end{corollary}

\section{Topological indices of vertex coalescence of graphs}\label{top}
In this section, some topological indices of vertex coalescence of graphs are computed. Using these we calculate the indices of some family of graphs. Also, as an application, we compute the indices of the organic compound 1,2-dicyclohexylethane(\ce{C_{14}H_{26}}).

\begin{proposition}
Wiener index of $G_1\circ_1G_2$ is 
$$W(G_1\circ_1G_2)=W(G_1)+W(G_2)+(n_2-1)d_{G_1}(v)+(n_1-1)d_{G_2}(v),$$
where $v$ is the vertex that is merged in $G_1\circ_1G_2$.
\end{proposition}

\begin{proof}
Let $G=G_1\circ_1G_2$ and $v$ be the vertex merging in $G$. From Definition \ref{W},
\begin{align*}
W(G)=&\displaystyle\sum_{\{u,w\}\in V(G_1)}d(u,w)+\sum_{\{u,w\}\in V(G_2)}d(u,w)+\sum_{\substack{u\in V(G_1)\\w\in V(G_2)}}d(u,w)\\
=&W(G_1)+W(G_2)+(n_2-1)d_{G_1}(v)+(n_1-1)d_{G_2}(v). 
\end{align*}
\end{proof}
\begin{remark}
    The Wiener index of cycle is $W(C_m)=\begin{cases}
        \frac{m^3}{8} & \text{if $m$ is even}\\
        \frac{m(m^2-1)}{8} & \text{if $m$ is odd,}
    \end{cases}$ and Wiener index of path is $W(P_n)=\frac{n(n^2-1)}{6}$. Thus the Wiener index of Lollipop graph $L(m,n-1)$ is  $$W(L(m,n-1))=\begin{cases}
        \frac{m^3}{8}+\frac{n(n^2-1)}{6}+(n-1)\left(\frac{m^2+2n(m-1)}{4}\right) & \text{if $m$ is even}\\
        \frac{m(m^2-1)}{8}+\frac{n(n^2-1)}{6}+\frac{(n-1)(m-1)(m+1+2n)}{4} & \text{if $m$ is odd.}
    \end{cases}$$
\end{remark}

    \begin{remark}
    The Dumbbell graph, denoted by $D_{l,m,n-3}$, is obtained from the coalescence of a cycle $C_l$ and the pendant vertex of a Lollipop graph $L(m,n-1)$.
    
    The Wiener index of Dumbbell graph $D_{m,m,n-3}$ is 
     $$W(D_{m,m,n-3})=\begin{cases}
        \frac{m^3}{4}+\frac{n(n^2-1)}{6}+\frac{m(m^2+3mn-4m+4)+n(4-6m+2mn-2n)-2}{2} & \text{if $m$ is even}\\
        \frac{m(m^2-1)}{4}+\frac{n(n^2-1)}{6}+(m-1)\frac{m^2-3m+3mn-3n+4n^2}{2} & \text{if $m$ is odd.}
    \end{cases}$$
\end{remark}

\begin{figure}[H]
    \centering
  \begin{tikzpicture}  
  [scale=0.9,auto=center]
   \tikzset{dark/.style={circle,fill=black}}

 \node[dark] (d1) at (0,1.5){} ;  
  \node[dark] (d2) at (1.5,0)  {}; 
  \node[dark] (d3) at (3,1.5)  {};  
  \node[dark] (d4) at (1.5,3) {};  
   \node[dark] (d5) at (4.5,1.5)  {}; 
   \node[dark] (d6) at (6,1.5)  {}; 
   \node[dark] (d7) at (7.5,1.5)  {}; 

   \node[dark] (d8) at (8.5,3){} ;  
  \node[dark] (d9) at (8.5,0)  {}; 
  
  \node[dark] (d10) at (10,3)  {};  
  \node[dark] (d11) at (10,0) {};  
  
   \node[dark] (d12) at (11,1.5)  {};

   \draw (d1) -- (d2);
  \draw (d2) -- (d3);  
  \draw (d3) -- (d4);  
  \draw (d4) -- (d1);  
  \draw (d3) -- (d5);
  \draw (d5) -- (d6);
  \draw (d6) -- (d7);
   \draw (d7) -- (d8);
  \draw (d8) -- (d10);  
  \draw (d10) -- (d12);  
  \draw (d12) -- (d11);  
  \draw (d11) -- (d9);
  \draw (d9) -- (d7);


\end{tikzpicture}  
    \caption{$D_{4,6,1}$}
    \label{dumbbell}
\end{figure}

\begin{example}
    The molecular graph of 1,2-dicyclohexylethane(\ce{C_{14}H_{26}}) is isomorphic to $D_{6,6,1}$. Therefore the Wiener index of 1,2-dicyclohexylethane(\ce{C_{14}H_{26}}) is 343.
    \begin{figure}[H]
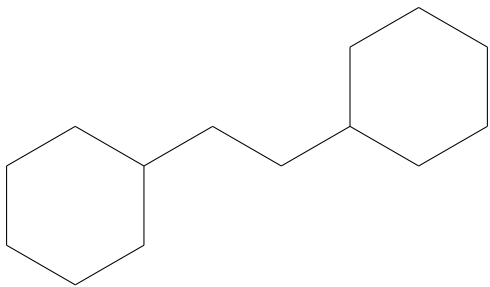

        \centering
         \chemfig{*6(---((-[:30]-[:-30]-[:30](*6(------))))---)}
        \caption{1,2-dicyclohexylethane(\ce{C_{14}H_{26}})}
        \label{dicyclohexylethane}
    \end{figure}
   
\end{example}
 \begin{remark}
    The Dandelion graph, denoted by $D(m,n)$, is obtained from the coalescence of a Star $S_n$ and the pendant vertex of a path $P_m$.
    
    The Wiener index of Dandelion graph $D(m,n)$ is 
     $$W(D(m,n))=(n-1)^2+\frac{m(m^2-1)}{6}+\frac{(n-1)(m-1)(m+2)}{2}.$$
\end{remark}

\begin{figure}[htbp]
\centering
\begin{tikzpicture}
   [scale=0.5,auto=center]
   \tikzset{dark/.style={circle,fill=black, scale=0.03cm}}

    \node[dark] (d2) at (1.5,3.7){};
      \node[dark] (d3) at (3.7,1.6){};
       \node[dark] (d4) at (3.7,-1.5){};
        \node[dark] (d5) at (1.5,-3.7){};
         \node[dark] (d6) at (-1.5,-3.7){};
          \node[dark] (d7) at (-3.7,-1.4){};
           \node[dark] (d8) at (-3.8,1.4){};
            \node[dark] (d9) at (-1.6,3.7){};
  \node[dark] (d10) at (0,4){};
   \node[dark] (d11) at (0,-4){};
    \node[dark] (d12) at (4,0){};
     \node[dark] (d13) at (-4,0){};
     \node[dark] (d14) at (0,0){};   

               \node[dark] (d15) at (2.8,2.8){};  
     \node[dark] (d16) at (-2.8,-2.8){};  
        \node[dark] (d17) at (2.8,-2.8){};
         \node[dark] (d18) at (-2.8,2.8){};
          \node[dark] (d19) at (8,0){};
           \node[dark] (d20) at (12,0){};
            \node[dark] (d21) at (16,0){};

            \draw (d14) -- (d2);
             \draw (d14) -- (d3);
              \draw (d14) -- (d4);
                \draw (d14) -- (d5);
             \draw (d14) -- (d6);
              \draw (d14) -- (d7);
                \draw (d14) -- (d8);
             \draw (d14) -- (d9);
              \draw (d14) -- (d10);
                \draw (d14) -- (d11);
             \draw (d14) -- (d12);
              \draw (d14) -- (d13);  
             \draw (d14) -- (d15);
              \draw (d14) -- (d16);
                \draw (d14) -- (d17);
             \draw (d14) -- (d18);
              \draw (d12) -- (d19);
        \draw (d19) -- (d20);
        \draw (d21) -- (d20);

\end{tikzpicture}
\caption{$D(16,5)$.}
 \end{figure}

  \begin{remark}
    The Kite graph, denoted by $Ki_{n,m}$, is obtained from the coalescence of a complete graph $K_n$ and the pendant vertex of a path $P_m$.
    
    The Wiener index of Kite graph $Ki_{n,m}$ is 
     $$W(Ki_{n,m})=\frac{n(n-1)}{2}+\frac{m(m^2-1)}{6}+\frac{(n-1)(m-1)(m+2)}{2}.$$
\end{remark}

\begin{figure}[htbp]
\centering
\begin{tikzpicture}
   [scale=0.5,auto=center]
   \tikzset{dark/.style={circle,fill=black, scale=0.03cm}}

    \node[dark] (d2) at (30,2.01){};
      \node[dark] (d3) at (28.09,0.62){};
       \node[dark] (d4) at (28.82,-1.63){};
        \node[dark] (d5) at (31.18,-1.63){};
         \node[dark] (d6) at (31.91,0.62){};
          \node[dark] (d7) at(34.41,0.62){};
           \node[dark] (d8) at (36.91,0.62){};
            \node[dark] (d9) at (39.41,0.62){};

            \draw (d3) -- (d2);
             \draw (d4) -- (d3);
              \draw (d5) -- (d4);
                \draw (d6) -- (d5);
             \draw (d2) -- (d4);
              \draw (d2) -- (d5);
                \draw (d2) -- (d6);
             \draw (d3) -- (d5);
              \draw (d3) -- (d6);
                \draw (d4) -- (d6);
             \draw (d6) -- (d7);
              \draw (d7) -- (d8);  
             \draw (d8) -- (d9);

\end{tikzpicture}
\caption{$Ki_{5,4}$.}
 \end{figure}

\begin{proposition}
Hyper-Wiener index of $G_1\circ_1G_2$ is 
\begin{align*}
WW(G_1\circ_1G_2)&=WW(G_1)+WW(G_2)\\
&+\frac{1}{2}\Big((n_2-1)(d_{G_1}(v)+d^2_{G_1}(v))+(n_1-1)(d_{G_2}(v)+d^2_{G_2}(v))+2d_{G_1}(v)d_{G_2}(v)\Big),
\end{align*}
where $v$ is the vertex that is merged in $G_1\circ_1G_2$.
\end{proposition}

\begin{proof}
Let $G=G_1\circ_1G_2$ and $v$ be the vertex merging in $G$. From Definition \ref{WW},
\begin{align*}
    WW(G)=&\frac{1}{2}(W(G_1)+W(G_2)+(n_2-1)d_{G_1}(v)+(n_1-1)d_{G_2}(v))\\
    &+\frac{1}{2}\left(\sum_{\{u,w\}\in V(G_1)}d^2(u,w)+\sum_{\{u,w\}\in V(G_2)}d^2(u,w)+\sum_{\substack{u\in V(G_1)\\w\in V(G_2)}}d^2(u,w)\right)\\
    =&WW(G_1)+WW(G_2)\\
&+\frac{1}{2}\Big((n_2-1)(d_{G_1}(v)+d^2_{G_1}(v))+(n_1-1)(d_{G_2}(v)+d^2_{G_2}(v))+2d_{G_1}(v)d_{G_2}(v)\Big).
\end{align*}
\end{proof}

\begin{remark}
    The hyper-Wiener index of of Lollipop graph $L(m,n-1)$ is\\
    
    $WW(L(m,n-1))=$
    $$\begin{cases}
        \frac{m^2(m+1)(m+2)}{48}+\frac{n^4+2n^3-n^2-2n}{24}+\frac{(n-1)(m(m^2+3m+2)+4n(m-1)(n+1)+3m^2n)}{24} & \text{if $m$ is even}\\
        \frac{m(m^2-1)(m+3)}{48}+\frac{n^4+2n^3-n^2-2n}{24}+\frac{(m-1)(n-1)((m+1)(m+3)+4n(n+1)+3n(m+1))}{24} & \text{if $m$ is odd.}
    \end{cases}$$
    \end{remark}

    \begin{remark}
    The hyper-Wiener index of Dumbbell graph $D_{m,m,n-3}$ is \\

    $WW(D_{m,m,n-3})=$
     $$\begin{cases}
        \frac{m^2(m+1)(m+2)}{24}+\frac{n^4+2n^3-n^2-2n}{24}\\
        \hspace{1cm}+\frac{7 m^4+4 m^3(-5+7 n)-8 n\left(1-3 n+2 n^2\right)+4 m^2\left(2-12 n+9 n^2\right)+8 m\left(-2+5 n-6 n^2+2 n^3\right)}{48} & \text{if $m$ is even}\\
        \frac{m(m^2-1)(m+3)}{24}+\frac{n^4+2n^3-n^2-2n}{24}\\
        \hspace{1cm}+\frac{(-1+m)\left(-3+7 m^3-16 n-12 n^2+16 n^3+m^2(-13+28 n)+m\left(-23-20 n+36 n^2\right)\right)}{48}& \text{if $m$ is odd.}
    \end{cases}$$
\end{remark}

\begin{example}
    The hyper-Wiener index of 1,2-dicyclohexylethane is 1032.
\end{example}

 \begin{remark}
    The hyper-Wiener index of Dandelion graph $D(m,n)$ is 
     $$WW(D(m,n))=\frac{(n-1)(3n-4)}{2}+\frac{m(m+2)(m^2-1)}{24}+\frac{(n-1)(m-1)(m^2+4m+6)}{6}.$$
\end{remark}

 \begin{remark}
    The hyper-Wiener index of Kite graph $Ki_{m,n}$ is 
     $$WW(Ki_{m,n})=\frac{(n-1)(n-1)}{2}+\frac{m(m+2)(m^2-1)}{24}+\frac{(n-1)(m-1)(m^2+4m+6)}{6}.$$
\end{remark}

\begin{proposition}
The forgotten topological index of $G_1\circ_1G_2$ is 
$$F(G_1\circ_1G_2)=F(G_1)+F(G_2)+3deg_{G_1}(v)deg_{G_2}(v)(deg_{G_1}(v)+deg_{G_2}(v)),$$
where $v$ is the vertex that is merged in $G_1\circ_1G_2$.
\end{proposition}

\begin{proof}
Let $G=G_1\circ_1G_2$ and $v$ be the vertex merging in $G$.
From Definition \ref{F},
\begin{align*}
F(G)=&\sum_{u\in V(G_1)}deg_{G_1}^3(u)-deg_{G_1}^3(v)+\sum_{u\in V(G_2)}deg_{G_2}^3(u)-deg_{G_2}^3(v)\\
&\hspace{2cm}+\left(deg_{G_1}(v)+deg_{G_2}(v)\right)^3\\
=&F(G_1)+F(G_2)+3deg_{G_1}(v)deg_{G_2}(v)(deg_{G_1}(v)+deg_{G_2}(v)).
\end{align*}
\end{proof}

\begin{remark}
    The forgotten topological index of of Lollipop graph $L(m,n-1)$ is  $$F(L(m,n-1))=8(m+n)+4.$$
    \end{remark}

    \begin{remark}
    The forgotten topological index of Dumbbell graph $D_{m,m,n-3}$ is 
     $$F(D_{m,m,n-3})=8(2m+n)+22.$$
\end{remark}

\begin{example}
    The forgotten topological index of 1,2-dicyclohexylethane is 150.
\end{example}

 \begin{remark}
    The forgotten topological index of Dandelion graph $D(m,n)$ is 
     $$F(D(m,n))=n^3+n-16+8m.$$
\end{remark}

 \begin{remark}
    The forgotten topological index of Kite graph $Ki_{m,n}$ is 
     $$F(Ki_{m,n})=n(n-1)((n-1)^2+3)-14+8m.$$
\end{remark}

\begin{proposition}
First Zagreb index of $G_1\circ_1G_2$ is 
$$M_1(G_1\circ_1G_2)=M_1(G_1)+M_1(G_2)+2deg_{G_1}(v)deg_{G_2}(v),$$
where $v$ is the vertex that is merged in $G_1\circ_1G_2$.
\end{proposition}

\begin{proof}
Let $G=G_1\circ_1G_2$ and $v$ be the vertex merging in $G$.
From Definition \ref{M1},
\begin{align*}
M_1(G)&=\sum_{u\in V(G_1)}deg_{G_1}^2(u)-deg_{G_1}^2(v)+\sum_{u\in V(G_2)}deg_{G_2}^2(u)-deg_{G_2}^2(v)\\
&\hspace{2cm}+\left(deg_{G_1}(v)+deg_{G_2}(v)\right)^2\\
&=M_1(G_1)+M_1(G_2)+2deg_{G_1}(v)deg_{G_2}(v).
\end{align*}
\end{proof}

\begin{remark}
    The first Zagreb index of of Lollipop graph $L(m,n-1)$ is  $$M_1(L(m,n-1))=4(m+n)-2.$$
    \end{remark}

    \begin{remark}
    The first Zagreb index of Dumbbell graph $D_{m,m,n-3}$ is 
     $$M_1(D_{m,m,n-3})=4(2m+n)+2.$$
\end{remark}

\begin{example}
    The first Zagreb index of 1,2-dicyclohexylethane is 66.
\end{example}

  \begin{remark}
    The first Zagreb index of Dandelion graph $D(m,n)$ is 
     $$M_1(D(m,n))=n^2+n-8+4m.$$
\end{remark}

 \begin{remark}
    The first Zagreb index of Kite graph $Ki_{m,n}$ is 
     $$M_1(Ki_{m,n})=(n-1)(n^2-n+2)+4m-6.$$
\end{remark}

\begin{proposition}
Narumi-Katayama index of $G_1\circ_1G_2$ is 
$$NK(G_1\circ_1G_2)=NK(G_1)NK(G_2)\frac{deg_{G_1}(v)+deg_{G_2}(v)}{deg_{G_1}(v)deg_{G_2}(v)},$$
where $v$ is the vertex that is merged in $G_1\circ_1G_2$.
\end{proposition}

\begin{proof}
Let $G=G_1\circ_1G_2$ and $v$ be the vertex merging in $G$.
From Definition \ref{NK},
\begin{align*}
NK(G)&=\frac{\displaystyle\prod_{u\in V(G_1)}deg_{G_1}(u)\prod_{u\in V(G_2)}deg_{G_2}(u)}{deg_{G_1}(v)deg_{G_2}(v)}\left(deg_{G_1}(v)+deg_{G_2}(v)\right)\\
&=NK(G_1)NK(G_2)\frac{deg_{G_1}(v)+deg_{G_2}(v)}{deg_{G_1}(v)deg_{G_2}(v)}.
\end{align*}
\end{proof}

\begin{remark}
    The Narumi-Katayama index of of Lollipop graph $L(m,n-1)$ is  $$NK(L(m,n-1))=3\times2^{m+n-3}.$$
    \end{remark}

    \begin{remark}
    The Narumi-Katayama index of Dumbbell graph $D_{m,m,n-3}$ is 
     $$NK(D_{m,m,n-3})=9\times 2^{2m+n-4}.$$
\end{remark}

\begin{example}
    The Narumi-Katayama index of 1,2-dicyclohexylethane is 36864.
\end{example}
\begin{remark}
    The Narumi-Katayama index of Dandelion graph $D(m,n)$ is 
     $$NK(D(m,n))=n.2^{m-2}.$$
\end{remark}

\begin{remark}
    The Narumi-Katayama index of Kite graph $Ki_{m,n}$ is 
     $$NK(Ki_{m,n})=n(n-1)^{n-1}2^{m-2}.$$
\end{remark}



\section{Conclusion}
This paper estimates some structural properties of a non-regular graph obtained from the $k$-coalescence of two graphs. Also, the $A_\alpha$-characteristic polynomial of $k$-coalescence of two graphs is determined. Moreover, the $A_\alpha$-spectrum and $A_\alpha$-energy of $k$-coalescence of two complete graphs are computed. In addition, some topological indices of vertex coalescence of two graphs are estimated. The Wiener index, hyper-Wiener ined etc. of some family of graphs are derived as an application. Using these, the indices of the organic compound 1,2-dicyclohexylethane(\ce{C_{14}H_{26}}) are calculated.

\bibliographystyle{unsrt}  
\bibliography{references}

\end{document}